\documentclass[12pt]{amsart}
\usepackage[all]{xypic}
\usepackage{amsmath,graphics}
\usepackage{amsfonts,amssymb}

\theoremstyle{plain}
\newtheorem*{theorem*}{Theorem}
\newtheorem*{lemma*} {Lemma}
\newtheorem*{corollary*} {Corollary}
\newtheorem*{proposition*} {Proposition}
\newtheorem{theorem}{Theorem}[section]
\newtheorem{lemma}[theorem]{Lemma}
\newtheorem{corollary}[theorem]{Corollary}
\newtheorem{proposition}[theorem]{Proposition}

\theoremstyle{remark}
\newtheorem*{remark}{Remark}
\newtheorem*{definition}{Definition}

\newtheorem*{claim}{Claim}

\theoremstyle{definition}

\def \K {\mathbf{K}}
\def \Z {\mathbf{Z}}

\def \F {\mathbf{F}}

\def\R{\Bbb{R}}

\def\genus{\mbox{genus}}
\def\eps{\epsilon}

\def\s{\sigma}

\def\gl{\mbox{GL}}

\def\k{\Bbb{K}}
\def\K{\Bbb{K}}

\def\id{\mbox{id}}

\def\Z{\Bbb{Z}}

\def\N{\Bbb{N}}
\def\l{\lambda}
\def\part{\partial}

\def\U{\mathcal{U}}

\def\a{\alpha}
\def\g{\gamma}

\def\bp{\begin{pmatrix}}

\def\sm{\setminus}
\def\ep{\end{pmatrix}}
\def\bn{\begin{enumerate}}

\def\rank{\mbox{rank}}

\def\en{\end{enumerate}}
\def\ba{\begin{array}}
\def\ea{\end{array}}

\def\s{\sigma}
\def\a{\alpha}

\def\ti{\tilde}

\def\fr12{\frac{1}{2}}

\def\im{\mbox{Im}}

\def\ker{\mbox{Ker}}

\def\ord{\mbox{ord}}

\def\hom{\mbox{Hom}}

\def\kt{\K_{\g}\tpm}
\def\dimm{\mbox{dim}}
\def\deg{\mbox{deg}}

\def\gcd{\mbox{gcd}}

\def\tpm{[t^{\pm 1}]}

\def\gen{\mbox{genus}}

\def\zt{\Z\tpm}

\def\f{\Bbb{F}}
\def\F{\Bbb{F}}
\def\k{\Bbb{K}}
\def\K{\Bbb{K}}
\def\kk{\Bbb{K}^d}
\def\kat{\k_{\g}[t^{\pm 1}]}
\def\ktfield{\k_{\g}(t)}
\def\kkat{\k_{\g}[t^{\pm 1}]^d}
\def\glkatk{\mbox{GL}(\k_{\g}[t^{\pm 1}],d)}

\def\ol{\overline}

\def\ft{\F[t^{\pm 1}]}
\def\fd{\F^d}

\def\zt{\Z[t^{\pm 1}]}

\def\K{\Bbb{K}}

\def\tnphi{||\phi||_T}

\def\gen{\mbox{genus}}

\def\degtaum{\deg(\tau(M,\a))}

\def\cmtbf#1{} \def\cmt#1{}
\def\glkatk{\mbox{GL}(\k_{\g}[t^{\pm 1}],d)}
\def\glkk{\mbox{GL}(\k,d)}

\begin{document}

\title{Reidemeister torsion, the Thurston norm and Harvey's invariants} \author{Stefan Friedl}
\date{\today} \address{Rice University, Houston, Texas, 77005-1892} \email{friedl@math.rice.edu}
\def\subjclassname{\textup{2000} Mathematics Subject Classification}
\expandafter\let\csname subjclassname@1991\endcsname=\subjclassname \expandafter\let\csname
subjclassname@2000\endcsname=\subjclassname 

\date{\today}
\begin{abstract}
Cochran introduced Alexander polynomials over non--commutative Laurent polynomial rings. Their
degrees were studied by Cochran, Harvey and Turaev as they give lower bounds on the Thurston
norm. We first extend Cochran's definition to twisted Alexander polynomials. We then show how
Reidemeister torsion relates to these invariants and we give lower bounds on the Thurston norm
in terms of the Reidemeister torsion. This gives in particular  a concise formulation of the
bounds of Cochran, Harvey and Turaev. The Reidemeister torsion approach also gives a natural
approach to proving and extending certain monotonicity results of Cochran and Harvey.
 \end{abstract}
\maketitle

\section{Introduction}
The following algebraic setup allows us to define twisted non--commutative Alexander
polynomials. First let $\k$ be a (skew) field and $\g:\k\to \k$ a ring homomorphism. Then denote
by $\kat$ the skew Laurent polynomial ring over $\k$. More precisely the elements in $\kat$ are
formal sums $\sum_{i=m}^n a_it^i$ ($m\leq n\in \Z$) with $a_i\in \k$. Addition is given by
addition of the coefficients, and multiplication is defined using the rule $t^ia=\g^i(a)t^i$ for
any $a\in \k$.

Let $X$ be a connected CW--complex with finitely many cells in dimension $i$. Given a
representation $\pi_1(X)\to \gl(\kt,d)$ we can consider the $\kat$--modules $H_i^\a(X;\kkat)$
and we define twisted non--commutative Alexander polynomials $\Delta_i^\a(t) \in \kat$ (cf.
Section \ref{section:alex} for details). Twisted Alexander polynomials over commutative Laurent
polynomial rings were first introduced by Lin \cite{Lin01}, Alexander polynomials over skew
Laurent polynomial rings were introduced by Cochran \cite{Co04}. Our definition is a combination
of the definitions in \cite{KL99} and \cite{Co04}. In Theorem \ref{thm:welldef} we describe the
indeterminacy of these polynomials.

Denote by $\ktfield$ the quotient field of $\kat$.
 We denote the induced representation $\pi_1(X)\to \gl(\kt,d)\to \gl(\ktfield,d)$ by $\a$
as well. If the homology groups $H_i^\a(X;\ktfield^d)$ vanish and if $X$ is a finite connected
CW--complex, then we can define the Reidemeister torsion $\tau(X,\a) \in K_1(\ktfield)/\pm
\a(\pi_1(X))$. An important tool is the Dieudonn\'e determinant which defines an isomorphism
\[  \det: K_1(\ktfield)\to \ktfield^\times_{ab},\]
where $\ktfield^\times_{ab}$ denotes the abelianization of the multiplicative group
$\ktfield^\times=\ktfield\sm \{0\}$. We can therefore study $\det(\tau(X,\a))\in
\ktfield^\times_{ab}/\pm \det(\a(\pi_1(X)))$. We refer to Sections \ref{section:torsion} and
\ref{section:dieudonne} for details.
 The following result generalizes well--known commutative results of Turaev
\cite{Tu86,Tu01} and Kirk--Livingston \cite{KL99}.

 \begin{theorem} \label{thm:taualexgeneral}
Let $X$ be a finite connected CW complex of dimension $n$. Let $\a:\pi_1(X)\to \gl(\kt,d)$ be a
 representation such that $H_*^\a(X;\ktfield^d)=0$. Then $\Delta_i^\a(t)\ne 0$ for all $i$ and
\[ \det(\tau(X,\a))=\prod_{i=0}^{n-1}\Delta_i^\a(t)^{(-1)^{i+1}} \in  \ktfield^\times_{ab}/\{kt^e|k\in \K\sm \{0\},e\in \Z\}.\]
\end{theorem}


For $f(t)=\sum_{i=m}^n a_it^i\in \kat\sm \{0\}$ with $a_m\ne 0, a_n \ne 0$, we define its degree
to be $\deg(f(t))=n-m$. We can extend this to a degree function $\deg:\ktfield\sm \{0\}\to \Z$.
We denote $\deg(\det(\tau(X,\a)))$ by $\deg(\tau(X,\a))$. Theorem \ref{thm:taualexgeneral} then
implies that the degree of $\tau(X,\a)$ is the alternating sum of the degrees of the twisted
Alexander polynomials (cf. Corollary \ref{cor:taudelta}).

We now turn to the study of 3--manifolds. Here and  throughout the paper we will assume that all
manifolds are compact, orientable and connected. Recall that given a  3--manifold $M$ and
$\phi\in H^1(M;\Z)$ the \emph{Thurston norm} (\cite{Th86}) of $\phi$ is defined as
 \[
||\phi||_{T}=\min \{ -\chi(\hat{S})\, | \, S \subset M \mbox{ properly embedded surface dual to
}\phi\}
\] where $\hat{S}$ denotes the result of discarding all connected components of $S$ with positive Euler
characteristic.
  As an example
consider $X(K)=S^3\sm \nu K$, where $K\subset S^3$ is a knot and $\nu K$ denotes an open tubular
neighborhood of $K$ in $S^3$. Let $\phi\in H^1(X(K);\Z)$ be a generator, then $||\phi||_T=2\,
\gen(K)-1$.

 Let $X$ be a connected CW--complex and let $\phi
\in H^1(X;\Z)$. We identify henceforth $H^1(X;\Z)$ with $\hom(H_1(X;\Z),\Z)$ and
$\hom(\pi_1(X),\Z)$. A representation $\a:\pi_1(X) \to \gl(\kat,d)$ is called
\emph{$\phi$--compatible} if for any $g\in \pi_1(X)$ we have $\a(g)=At^{\phi(g)}$ for some $A\in
\gl(\k,d)$. This generalizes a notion of Turaev \cite{Tu02b}.

The following theorem gives lower bounds on the Thurston norm using Reidemeister torsion. It
contains the lower bounds  of McMullen \cite{Mc02}, Cochran \cite{Co04}, Harvey \cite{Ha05},
Turaev \cite{Tu02b} and of the author together with Taehee Kim \cite{FK05}. To our knowledge
this theorem is the strongest of its kind. Not only does it contain these results, the
formulation of the inequalities in \cite{Mc02,Co04,Ha05,Tu02b} in terms of the degrees of
Reidemeister torsion also gives a very concise reformulation of their results.

\begin{theorem}\label{mainthm1}
Let $M$ be a 3--manifold with empty or toroidal boundary. Let $\phi \in H^1(M;\Z)$ and
$\a:\pi_1(M)\to \gl(\kat,d)$ a $\phi$--compatible representation. Then $\tau(M,\a)$ is defined
if and only if $\Delta^\a_1(t)\ne 0$. Furthermore if $\tau(M,\a)$ is defined, then
\[ ||\phi||_T \geq \frac{1}{d}\degtaum.\]
If $(M,\phi)$ fibers over $S^1$, then
\[ ||\phi||_T = \max\{0,\frac{1}{d}\degtaum\}.\]
\end{theorem}

The most commonly used skew fields are the quotient fields $\K(G)$ of group rings $\F[G]$ ($\F$
a commutative field) for certain torsion--free groups $G$, we refer to Section
\ref{section:skewfield} for details. The following theorem says roughly that `larger groups give
better bounds on the Thurston norm'. The main idea of the proof is to use the fact that
Reidemeister torsion behaves well under ring homomorphisms, in contrast to Alexander
polynomials. We refer to Section \ref{section:comparison} or to \cite{Ha06} for the definition
of an admissible triple.

\begin{theorem}\label{thm:monoton} \label{mainthm2}
Let $M$ be a 3--manifold with empty or toroidal boundary or let $M$ be a 2--complex with
$\chi(M)=0$. Let $\phi\in H^1(M;\Z)$. Let $\a:\pi_1(M)\to \gl(\F,d)$, $\F$ a commutative field,
be a representation and $(\varphi_G:\pi\to G,\varphi_H:\pi \to H,\phi)$ an admissible triple for
$\pi_1(M)$, in particular we have epimorphisms $G\to H\to \Z$. Write $G'=\ker\{G\to \Z\}$ and
$H'=\ker\{H\to \Z\}$.

If $\tau(M,\varphi_H\otimes \a)\in K_1(\K(H')(t))$ is defined, then $\tau(M,\varphi_G\otimes
\a)\in K_1(\K(G')(t))$ is defined. Furthermore in that case
\[ \deg(\tau(M,\varphi_G\otimes \a))\geq  \deg(\tau(M,\varphi_H\otimes \a)).\]
\end{theorem}

A similar theorem holds for 2--complexes with Euler characteristic zero. As a special case
consider the case that $\a$ is the trivial representation. Using Theorem
\ref{thm:taualexgeneral} we can recover the monotonicity results of \cite{Co04} and \cite{Ha06}.
We hope that our alternative proof using Reidemeister torsion will contribute to the
understanding of their results.
\\

The paper is organized as follows. In Section \ref{section:rt} we recall the definition of
Reidemeister torsion. In Section
 \ref{section:rtaltex} we introduce twisted non--commutative Alexander polynomials, we compute their
 indeterminacies
 in Theorem \ref{thm:welldef} and we
prove Theorem \ref{thm:taualexgeneral}. Beginning with Section \ref{section:3mfd} we concentrate
on 3--manifolds. In particular in Section \ref{section:3mfd} we give the proof of Theorem
\ref{mainthm1}. In Section \ref{section:ex} we give examples of $\phi$--compatible
representations. In Section \ref{section:comparison} we prove Theorem \ref{mainthm2} and in
Section \ref{section:Harvey} we show that it implies Cochran's and Harvey's monotonicity
results. We conclude with a few open questions in Section \ref{section:questions}.
\\

{\bf Acknowledgment:} The author would like to thank Stefano Vidussi for pointing out to him the
functoriality of  Reidemeister torsion and he would like to thank Tim Cochran, Shelly Harvey and
Taehee Kim for helpful discussions. The author also would like to thank the referee for many
helpful comments.


\section{Reidemeister torsion} \label{section:rt}
\subsection{Definition of $K_1(R)$} \label{section:defk1}
For the remainder of the paper we will only consider associative rings $R$ with $1\ne 0$ and
with the property that if $r\ne s \in \N_0$, then $R^r$ is not isomorphic to $R^s$ as an
$R$--module.

For such a ring $R$ define $\gl(R)=\underset{\rightarrow}{\lim} \, \gl(R,d)$, where we have the
following maps in the direct system: $\gl(R,d)\to \gl(R,d+1)$ given by $A \mapsto \bp A
&0\\0&1\ep $.
Then $K_1(R)$ is defined as $\gl(R)/[\gl(R),\gl(R)]$. Note that $K_1(R)$ is an abelian group.
For details we refer to \cite{Mi66} or \cite{Tu01}. There exists a canonical map $\gl(R,d)\to
K_1(R)$ for every $d$, in particular there exists a homomorphism from the units of $R$ into
$K_1(R)$. By abuse of notation we denote the image of $A\in \gl(R,d)$ in $K_1(R)$ by $A$ as
well. We denote by $-A$ the product of $A\in K_1(R)$ by the image of $(-1)$ under the map
$\gl(R,1)\to K_1(R)$.

We will often make use of the observation (cf. \cite[p.~61]{Ro94}) that for $A\in \gl(R,d_1),
B\in \gl(R,d_2)$ the product $AB\in K_1(R)$ is given by
\[  AB=\bp A &0\\0&B\ep \in K_1(R).\]



\subsection{Definition of Reidemeister torsion}
Let $C_*$ be a finite free  chain complex of $R$--modules. By this we mean a chain complex of
free finite right $R$--modules such that $C_i=0$ for all but finitely many $i\in \Z$. Let
$\mathcal{C}_i \subset C_i$ be a basis for all $i$ with $C_i\ne 0$. Assume that
$B_i=\im(C_{i+1})\subset C_{i}$ is free, pick a basis $\mathcal{B}_i$ of $B_i$ and a lift
$\ti{\mathcal{B}}_i$ of $\mathcal{B}_i$ to $C_{i+1}$. We write
$\mathcal{B}_i\ti{\mathcal{B}}_{i-1}$ for the collection of elements given by $\mathcal{B}_i$
and $\ti{\mathcal{B}}_{i-1}$. Since $C_*$ is acyclic this is indeed a basis for $C_i$. Then we
define the {\em Reidemeister torsion} of the based acyclic complex $(C_*,\{\mathcal{C}_i\})$ to
be
\[ \tau(C_*,\{\mathcal{C}_i\})=\prod  [\mathcal{B}_i\ti{\mathcal{B}}_{i-1}/\mathcal{C}_i]^{(-1)^{i+1}} \in K_1(R), \]
where $[d/e]$ denotes the matrix of a basis change, i.e. $[d/e]=(a_{ij})$ where $d_i=\sum_j
a_{ji}e_j$. Note that in contrast to \cite{Tu01} we view vectors as column vectors. This means
that our matrix is the transpose of the matrix in  \cite[p.~1]{Tu01}.

It is easy to see that $\tau(C_*,\{\mathcal{C}_i\})$ is independent of the choice of
$\{\mathcal{B}_i\}$ and of the choice of the lifts $\ti{\mathcal{B}}_i$.   If the $R$--modules
$B_i$ are not free, then one can show that they are stably free and a stable basis will then
make the definition work again. We refer to \cite[p.~369]{Mi66} or \cite[p.~13]{Tu01} for the
full details.


\subsection{Reidemeister torsion of a CW--complex} \label{section:torsion}

Let $X$ be a connected CW--complex. Denote the universal cover of $X$ by $\ti{X}$. We view
$C_*(\ti{X})$, the chain complex of the universal cover,  as a chain complex of right
$\Z[\pi_1(X)]$--modules, where the $\Z[\pi_1(X)]$--module structure is given via deck
transformations.

Let $R$ be a ring. Let $\a:\pi_1(X)\to \gl(R,d)$ be a representation. This equips $R^d$ with a
left $\Z[\pi_1(X)]$--module structure. We can therefore consider the  chain complex
$C_*^\a(X;R^d)=C_*(\ti{X})\otimes_{\Z[\pi_1(X)]} R^d$. Note that this is a finite free chain
complex  of (right) $R$--modules. We denote its homology by $H_*^\a(X;R^d)$, we drop $\a$ from
the notation if it is clear from the context.

Now assume that $X$ is a finite connected CW--complex. If $H_i^\a(X;R^d)=H_i(C_*^\a(X;R^d))\ne
0$ for some $i$ then we write $\tau(X,\a)=0$. Otherwise denote the $i$--cells of $X$ by
$\s_i^1,\dots,\s_{i}^{r_i}$ and denote by $e_1,\dots,e_d$ the standard basis of $R^d$. Pick an
orientation for each cell $\s_i^j$, and also pick a  lift $\ti{\s}_i^j$ for each cell $\s_i^j$
to the universal cover $\ti{X}$. We get a basis
 \[ \mathcal{C}_i=\{\ti{\s}^1_i\otimes e_1,\dots,\ti{\s}^1_i\otimes e_d,\dots,
 \ti{\s}_i^{r_i}\otimes e_1,\dots,\ti{\s}_i^{r_i}\otimes e_d\}\]
for $C_i^\a(X;R^d)$. Then we can define
\[ \tau(C_*^\a(X;R^d),\{\mathcal{C}_i\}) \in K_1(R).\]
This element depends only on the ordering and orientation of the cells and on the choice of
lifts of the cells to the universal cover. Therefore
\[ \tau(X,\a)=\tau(C_*^\a(X,R^d),\{\mathcal{C}_i\}) \in K_1(R)/\pm \a(\pi_1(X)) \]
is a well--defined invariant of the CW--complex $X$.

Now let $M$ be a compact PL--manifold. Pick any finite CW--structure for $M$ to define
$\tau(M,\a) \in K_1(R)/\pm \a(\pi_1(M))$. By Chapman's theorem \cite{Ch74} this is a
well--defined invariant of the manifold, i.e. independent of the choice of the CW--structure.

\subsection{Computation of Reidemeister torsion}\label{sec:comptau}
We explain an algorithm for computing Reidemeister torsion which was formulated by Turaev
\cite[Section~2.1]{Tu01} in the commutative case.

In the following assume that we have a finite free chain complex of $R$--modules
\[ 0\to C_m\xrightarrow{\partial_m} C_{m-1}\to \dots \to C_1\xrightarrow{\partial_1} C_0\to 0.\]
 Let $A_i=(a_{jk}^i)$ be the matrix representing $\partial_i$ corresponding to the given bases.
Note that in contrast to \cite[Section~2.1]{Tu01} we view the elements in $R^{rank(C_i)}$ as
column vectors.

Following Turaev \cite[p.~8]{Tu01} we define a \emph{matrix chain} for $C$ to be a collection of
sets $\xi=(\xi_0,\xi_1,\dots,\xi_m)$ where $\xi_i\subset \{1,2,\dots,\rank(C_i)\}$ so that
$\xi_0=\emptyset$. Given a matrix chain $\xi$ we define $A_i(\xi)$, $i=1,\dots,m$ to be the
matrix formed by the entries $a_{jk}^i$ with $j\not\in \xi_{i-1}$ and $k\in \xi_{i}$. Put
differently the matrix $(a_{jk}^i)_{jk}$ is given by considering only the $\xi_i$--columns of
$A_i$ and with the $\xi_{i-1}$--rows removed.

We say that a matrix chain $\xi$ is a $\tau$--chain if $A_1(\xi),\dots,A_{m}(\xi)$ are square
matrices. The following is the generalization of \cite[Theorem~2.2]{Tu01} to the
non--commutative setting. Turaev's proof can easily be generalized to this more general setting.

\begin{theorem} \label{thm:comptauturaev}
Let $\xi$ be a $\tau$--chain such that $A_i(\xi)$ is invertible for all odd $i$. Then $A_i(\xi)$
is invertible for all even $i$ if and only if $H_*(C)=0$. Furthermore if $H_*(C)=0$, then
\[ \tau(C)=\eps \prod_{i=1}^{m} A_i(\xi)^{(-1)^{i}} \in K_1(R)\]
for some $\eps \in \{\pm 1\}$.
\end{theorem}

This proposition is the reason why Reidemeister torsion behaves in general well under ring
homomorphisms.

\section{Reidemeister torsion and Alexander polynomials} \label{section:rtaltex}

\subsection{Laurent polynomial rings and the Dieudonn\'e determinant}  \label{section:dieudonne}
For the remainder of this paper let $\k$ be a (skew) field and let $\kat$ be a skew Laurent
polynomial ring. By \cite[Corollary~6.3]{DLMSY03} the ring $\kat$ has a classical quotient field
$\ktfield$ which is flat over $\kt$ (cf. \cite[p.~99]{Ra98}). In particular we can view $\kat$
as a subring of $\ktfield$ and any element in $\ktfield$ is of the form $f(t)g(t)^{-1}$ for some
$f(t)\in \kat$ and $g(t)\in \kat\sm \{0\}$. We refer to Theorem \ref{thm:tfa} for a related
result. Recall that we write $\ktfield^\times=\ktfield\sm \{0\}$.

In the following we mean by an elementary column (row) operation the addition of a right
multiple (left multiple) of one column (row) to a different column (row). Let $A$ be an
invertible $k\times k$--matrix over the skew field $\ktfield$. After elementary row operations
we can turn $A$ into a diagonal matrix $D=(d_{ij})$. Then the Dieudonn\'e determinant
$\det(A)\in \ktfield_{ab}^\times=\ktfield^\times/[\ktfield^\times,\ktfield^\times]$ is defined
to be $\prod_{i=1}^k d_{ii}$.  This is a well--defined map. Note that the Dieudonn\'e
determinant is invariant under elementary row and column operations. The Dieudonn\'e determinant
induces an isomorphism $\det:K_1(\ktfield)\to \ktfield_{ab}^\times$.  Using the last observation
in Section \ref{section:defk1} it is easy to see that $A=\det(A)\in K_1(\ktfield)$. We will
often make use of this equality. We refer to \cite[Theorem~2.2.5 and Corollary~2.2.6]{Ro94} for
more details.

In the introduction we defined $\deg:\kt\sm\{0\}\to \N$. This can be extended to a homomorphism
$\deg:\ktfield^\times \to \Z$ via $\deg(f(t)g(t)^{-1})=\deg(f(t))-\deg(g(t))$ for $f(t),g(t)\in
\kt\sm\{0\}$.  Clearly the degree map vanishes on $[\ktfield^\times,\ktfield^\times]$ and we get
an induced homomorphism $K_1(\ktfield)\to \ktfield^\times_{ab}\to \Z$ which we also denote by
$\deg$.

\subsection{Orders of $\kt$--modules}
 Let $H$ be a finitely generated right $\kt$--module. The ring $\kat$ is a principal
ideal domain (PID) since $\k$ is a skew field. We can therefore find an isomorphism
\[ H\cong  \bigoplus_{i=1}^l
\kat/p_i(t)\kt\] for  $p_i(t)\in \kat$ for $i=1,\dots,l$. Following \cite{Co04}  we define
$\ord(H)=\prod_{i=1}^lp_i(t) \in \kat$. This is called the \emph{order of $H$}.

Note that $\ord(H)\in \kat$ has a high degree of indeterminacy. For example writing the $p_i(t)$
in a different order will change $\ord(H)$. Furthermore we can change $p_i(t)$ by multiplication
by any element of the form $kt^e$ where $k\in \K^\times=\K\sm\{0\}$ and $e\in \Z$.
 The
following theorem can be viewed as saying that these are all possible indeterminacies.

\begin{theorem} \label{thm:welldef}
Let $H$ be a finitely generated right $\kt$--module. Then $\ord(H)=0$ if and only if $H$ is not
$\kat$--torsion. If $\ord(H)\ne 0$, then $\ord(H)\in \kt$ is well--defined considered as an
element in $\ktfield^\times_{ab}$ up to multiplication by an element of the form $kt^e, k\in
\K^\times$ and $e\in \Z$.
\end{theorem}

The first statement is clear. We postpone the proof of the second statement of the theorem to
Section \ref{section:proofwelldef}. We refer to \cite[p.~367]{Co04} for an alternative
discussion of the indeterminacy of $\ord(H)$. Note that the idea of considering $\ord(H)$ as an
element in $\ktfield^\times_{ab}$ is already present in \cite[p.~367]{Co04}.

It follows from Theorem \ref{thm:welldef} that $\deg(\ord(H))$ is well--defined. In fact we have
the following interpretation of $\ord(H)$.
\begin{lemma} \label{lem:finite}\cite[p.~368]{Co04}
Let $H$ be a finitely generated right $\kt$--torsion module. Then
\[ \deg(\ord(H))=\dim_\K(H).\]
\end{lemma}

Here we used that by \cite[Proposition~I.2.3]{St75} and \cite[p.~48]{Coh85} every right
$\k$--module $V$ is free and has a well-defined dimension $\dim_\K(V)$.

\begin{proof}
It is easy to see that for $f(t)\in \kt\sm\{0\}$ we have
\[ \deg(f(t))=\dimm_{\K}(\kt /f(t)\kt). \]
The lemma is now immediate.
\end{proof}

\subsection{Alexander polynomials}  \label{section:alex}
Let $X$ be a connected CW--complex with finitely many cells in dimension $i$.  Let
$\a:\pi_1(X)\to \gl(\kat,d)$ be a  representation. The right $\kat$--module $H_i(X;\kkat)$ is
called twisted (non--commutative) Alexander module. Similar modules were studied in \cite{Co04},
\cite{Ha05} and \cite{Tu02b}. Note that $H_i(X;\kkat)$ is a finitely generated $\kt$--module
since we assumed that $X$ has only finitely many cells in dimension $i$ and since $\kt$ is a
PID. We now define $\Delta^\a_{i}(t)=\ord(H_i(X;\kkat))\in \kt$, this is called the (twisted)
{\em $i$--th Alexander polynomial} of $(X,\a)$.

The degrees of these polynomials (corresponding to one--dimensional representations) have been
studied  recently  in various contexts (cf. \cite{Co04,Ha05,Ha06,Tu02b,LM05,FK05b,FH06}). We
hope that by determining   the indeterminacy of the Alexander polynomials (Theorem
\ref{thm:welldef}) more information can be extracted from the Alexander polynomials than just
the degrees.

\subsection{Proof of Theorem \ref{thm:welldef}} \label{section:proofwelldef}
We first point out that $\kt$ is a Euclidean ring with respect to the degree function. This
means that given $f(t),g(t)\in \kt\sm\{0\}$ we can find $a(t),r(t)\in \kt$ such that
$f(t)=g(t)a(t)+r(t)$ and such that either $r(t)=0$ or $\deg(r(t))<\deg(g(t))$.

Let $A$ be an $r\times s$--matrix over $\kt$ of rank $r$. Here and in the following the rank of
a matrix over $\kt$ will be understood as the rank of the matrix considered as a matrix over the
skew field $\ktfield$. Note that $\rank(A)=r$ implies that in particular $s\geq r$. Since $\kt$
is a Euclidean ring we can perform a sequence of elementary row and column operations to
 turn $A$ into a matrix of the form $\bp D &0_{r\times (s-r)}\ep$ where $D$ is an $r\times r$--matrix
 and $0_{r\times (s-r)}$ stands for the $r\times (s-r)$--matrix consisting only of zeros.
 Since $A$ is of rank $r$ it follows that $D$ has rank $r$ as
well, in particular $D$ is a square matrix which is invertible over $\ktfield$ and we can consider
its Dieudonn\'e determinant $\det(D)$. We define $\det(A)=\det(D)\in \ktfield^\times_{ab}$.

\begin{lemma} \label{lem:detkt}
Let $A$ be a (square) matrix  over $\kt$ which is invertible over $\ktfield$.
\bn
\item
The Dieudonn\'e determinant $\det(A)\in \ktfield^\times_{ab}$ can be represented by an element
in $\kt\sm \{0\}$.
\item
If $A\in \gl(\kt,d)$, then $\det(A)\in \ktfield^\times_{ab}$ can be represented by an element of
the form  $kt^e, k\in \K^\times, e\in \Z$.
\item The Dieudonn\'e determinant induces a homomorphism
\[ \det:K_1(\ktfield) \to \ktfield^\times_{ab}\]
which sends $K_1(\kt)$ to $\{kt^e | k\in \K^\times,e\in
\Z\}/[\ktfield^\times,\ktfield^\times]\subset \ktfield^\times_{ab}$.
\en
\end{lemma}

\begin{proof}
The first statement follows from the discussion preceding the lemma. Now let $A\in \gl(\kt,r)$.
It follows from Lemma \ref{lem:finite} applied to $H=\kt^r/A\kt^r$ that $\deg(\det(A))=0$. This
proves the second statement. The last statement follows from the second statement and the fact
that the Dieudonn\'e determinant induces a
 homomorphism $\det:K_1(\ktfield) \to \ktfield^\times_{ab}$.
\end{proof}

\begin{proposition}
Let $A$ be an $r\times s$--matrix  over $\kt$ of rank $r$. Then $\det(A)\in
\ktfield^\times_{ab}$ is well--defined up to multiplication by an element of the form $kt^e,
k\in \k^\times,e\in \Z$. Furthermore $\det(A)$ is invariant under elementary row and column
operations.
\end{proposition}

\begin{proof}
First note that the effect of an elementary row operation on $A$  over $\kt$ can be described by
left multiplication by a matrix $P\in \gl(\kt,r)$. Similarly an elementary column operation on
$A$ over $\kt$ can be described by right multiplication by an $s\times s$--matrix $Q\in
\gl(\kt,s)$.

Now assume we have  $P_1,P_2\in\gl(\kt,r)$ and $Q_1,Q_2\in\gl(\kt,s)$ such that $P_iAQ_i=\bp
D_i&0_{r\times (s-r)}\ep, i=1,2$ where $D_i$ is an $r\times r$--matrix. We are done once we show
that $\det(D_1)=kt^e\det(D_2)\in \ktfield^\times_{ab}$ for some $ k\in \k^\times,e\in \Z$. Let
$E_i=P_i^{-1}D_i$. Then by Lemma \ref{lem:detkt} we only have to show that $E_1=E_2\in
K_1(\ktfield)/K_1(\kt)$.

We have $\bp E_1&0\ep Q_1^{-1}=\bp E_2&0\ep Q_2^{-1}$. Let $Q:=Q_2^{-1}Q_1\in \gl(\kt,s)$, we
therefore get
the equality $\bp E_1&0\ep=\bp E_2&0\ep Q$. Now write $Q=\bp Q_{11}&Q_{12}\\
Q_{21}&Q_{22}\ep$ where $Q_{ij}$ is a $n_i\times n_j$ matrix over $\kt$ with $n_1=r$ and
$n_2=s-r$. We get the equality
\[ \bp E_1 &0\\ Q_{21}&Q_{22}\ep = \bp E_2&0\\ 0&\id_{s-r}\ep \bp Q_{11}&Q_{12}\\
Q_{21}&Q_{22}\ep.\] It follows in particular that $Q_{22}$ is invertible over $\ktfield$.
Furthermore we have
\[ E_1 \cdot Q_{22}=E_2 \in K_1(\ktfield)/K_1(\kt).\]
 Note that $\deg:K_1(\ktfield)\to \Z$ vanishes on $K_1(\kt)$ by Lemma \ref{lem:detkt}.
 We therefore get $\deg(\det(E_1))+\deg(\det(Q_{22}))=\deg(\det(E_2))$, in
particular $\deg(\det(E_1))\leq \deg(\det(E_2))$. But by symmetry we have $\deg(\det(E_2))\leq
\deg(\det(E_1))$. In particular $\deg(\det(Q_{22}))=0$. The proposition now follows immediately
from  Lemma \ref{lem:detkt} since $\deg(f(t))=0$ for $f(t)\in \kt\sm \{0\}$ if and only if
$f(t)=kt^e$ for some $k\in \K^\times,e\in \Z$.

The last statement is immediate.
\end{proof}

Let $H$ be a finitely generated right $\kt$--module. We say that an $r\times s$--matrix $A$ is a
presentation matrix for $H$ if  the following sequence is exact:
\[  \kt^s\xrightarrow{A}\kt^r\to H\to 0.\]
We say that $A$ has full rank if the rank of $A$ equals $r$. Note that $A$ has full rank if and
only if $H\otimes_{\kt}\ktfield=0$.

The following lemma clearly implies Theorem \ref{thm:welldef}.

\begin{lemma} \label{lem:welldef}
Let $H$ be a finitely generated right $\kt$--module and let $A_1,A_2$ be presentation matrices
for $H$. Then $A_1$ has full rank if and only if $A_2$ has full rank. Furthermore if $A_i$ has
full rank, then
\[ \det(A_1)=\det(A_2)\in \ktfield^\times_{ab}/\{kt^e|k\in \K\sm\{0\}, e\in \Z\}. \]
\end{lemma}

\begin{proof}
It is well--known  that  any two presentation matrices for $H$ differ by a sequence of matrix moves
of the following forms and their inverses:
\bn
\item Permutation of rows or columns.
\item Replacement of the matrix $A$ by $\bp A&0\\ 0&1\ep$.
\item Addition of an extra column of zeros to the matrix $A$.
\item Addition of a right scalar multiple of a column to another column.
\item Addition of a left scalar multiple of a row to another row.
\en
This result is proven \cite[Theorem~6.1]{Li97} in the commutative case, but the proof carries
through in the case of the base ring $\kt$ as well (cf. also \cite[Lemma~9.2]{Ha05}).

 Clearly none of the moves changes the status of being of full
rank, and if a representation is of full rank, then it is easy to see that none of the moves
changes the determinant.
\end{proof}

\subsection{Proof of Theorem \ref{thm:taualexgeneral}} \label{section:proof11}



Now let $X$ be a finite connected CW complex of dimension $n$. Let $\a:\pi_1(X)\to \gl(\kt,d)$
be a
 representation such that $H_*^\a(X;\ktfield^d)=0$.
(Recall that we denote the induced representation $\pi_1(X)\to \gl(\ktfield,d)$ by $\a$ as
well).
 Furthermore recall that $\ktfield$
is flat over $\kt$, in particular $H_i(X;\ktfield^d)=H_i(X;\kt^d)\otimes_{\kt}\ktfield$. It
follows that $H_i(X;\ktfield^d)=0$ if and only if $H_i(X;\kkat)$ is $\kt$--torsion, which is
equivalent to $\Delta_{i}^{\a}(t)\ne 0$. This proves the first statement of Theorem
\ref{thm:taualexgeneral}. To conclude the proof of Theorem \ref{thm:taualexgeneral} it remains
to prove the following claim.

\begin{claim}
If $H_*^\a(X;\ktfield^d)=0$, then
\[ \det(\tau(X,\a))=\prod_{i=0}^{n-1}\Delta_i^\a(t)^{(-1)^{i+1}} \in  \ktfield^\times_{ab}/\{kt^e|k\in \K\sm \{0\},e\in \Z\}.\]
\end{claim}

\begin{proof}
Let $C_*=C_*(\ti{X})\otimes_{\Z[\pi_1(X)]}\kt^d$. Note that any $\kt$--basis for $C_*$ also
gives a basis for $C_*\otimes_{\kt}\ktfield$, which we will always denote by the same symbol.

Denote by  $\mathcal{C}_*$ the $\kt$--basis of $C_*$ as in Section \ref{section:torsion}.
 Let
$r_i:=\dim_{\K(t)}(C_i\otimes_{\kt} \ktfield)$ and let $s_i:=\dim_{\K(t)}(\ker\{C_i\otimes_{\kt}
\ktfield\xrightarrow{\partial_i} C_{i-1}\otimes_{\kt} \ktfield\})$. Note that
$s_{i}+s_{i-1}=r_i$ since $H_*(X;\ktfield^d)=0$. Furthermore note that
$\ker\{C_i\xrightarrow{\partial_i} C_{i-1}\}\subset C_i$ is a free direct summand of $C_i$ of
rank $s_i$ since $\kt$ is a PID. We can therefore
 pick  $\kt$--bases $\mathcal{C}'_i=\{v_1,\dots,v_{r_i}\}$ for $C_i$ such that
$\{v_1,\dots,v_{s_i}\}$ is a basis for $\ker\{C_i\xrightarrow{\partial_i} C_{i-1}\}$.  Note that
the base changes from $\mathcal{C}_i$ to $\mathcal{C}_i'$ are given by matrices which are
invertible over $\kt$. It follows that
\[ \tau(C_*\otimes_{\kt}\ktfield,\{\mathcal{C}_i\})=\tau(C_*\otimes_{\kt}\ktfield,\{\mathcal{C}'_i\}) \in K_1(\ktfield)/K_1(\kt).\]
Let $A_i$ be the $r_{i-1}\times r_i$--matrix representing $\partial_i:C_i\to C_{i-1}$ with
respect to the bases $\mathcal{C}'_i$ and $\mathcal{C}'_{i-1}$. Let
$\xi_i:=\{s_{i}+1,\dots,r_i\},i=1,\dots,n$ and $\xi_0:=\emptyset$. Let $A_i(\xi)$ as in Theorem
\ref{thm:comptauturaev}. Note that $A_i(\xi)$ is an $s_{i-1}\times s_{i-1}$--matrix over $\kt$.
In particular $\xi:=(\xi_0,\dots,\xi_n)$ is a $\tau$--chain. It is easy to see that
\[ A_i =\bp 0_{s_{i-1}\times s_i}&A_i(\xi)\\ 0_{s_{i-2}\times s_i}&0_{s_{i-2}\times s_{i-1}}\ep. \]
Since $A_i$ has rank $s_{i-1}$ it follows
 that $A_i(\xi)$ is  invertible over $\ktfield$. It
follows from Theorem \ref{thm:comptauturaev} that
\[ \tau(C_*\otimes_{\kt}\ktfield,\{\mathcal{C}'_i\}) =\prod_{i=1}^{n} A_i(\xi)^{(-1)^i}\in K_1(\ktfield).\]
We also have short exact sequences
\[ 0\to \kt^{s_{i-1}}\xrightarrow{A_{i}(\xi)} \kt^{s_{i-1}}\to H_{i-1}(C_*)=H_{i-1}(X;\kt^d)\to 0.\]
In particular $(A_i(\xi))$ is a presentation matrix for $H_{i-1}(X;\kt^d)$. It therefore follows
from Lemma \ref{lem:welldef} that $\det(A_i(\xi))=\Delta_{i-1}^\a(t)$.
\end{proof}

The following corollary now follows immediately from the fact that $\deg:\ktfield^\times\to \Z$
is a homomorphism and from Lemma \ref{lem:finite}.

\begin{corollary} \label{cor:taudelta}
Let $X$ be a finite connected CW complex of dimension $n$. Let $\a:\pi_1(X)\to \gl(\kt,d)$ be a
 representation such that $H_*^\a(X;\ktfield^d)=0$. Then
\[ \deg(\tau(X,\a))=\sum_{i=0}^{n-1}(-1)^{i+1}\deg(\Delta_i^\a(t))=\sum_{i=0}^{n-1}(-1)^{i+1}\dim_{\K}(H_i(X;\kt^d).\]
\end{corollary}

\begin{remark}
In the case that $H_*(X;\ktfield^d)\ne 0$ we can pick $\kt$--bases $\mathcal{H}_i$ for the
$\kt$--free parts of $H_i(X;\kt^d)$. These give  bases for
$H_i(X;\ktfield^d)=H_i(X;\kt^d)\otimes_{\kt}\ktfield$ and we can consider
$\tau(X,\a,\{\mathcal{H}_i\})=\tau(C_*(\ti{X})\otimes_{\Z[\pi_1(X)]}\ktfield^d,\{\mathcal{H}_i\})\in
K_1(\ktfield)$ (cf. \cite{Mi66} for details). If we consider $\tau(X,\a,\{\mathcal{H}_i\})\in
K_1(\ktfield)/K_1(\kt)$ then this is independent of the choice of $\{\mathcal{H}_i\}$. The proof
of Theorem \ref{thm:taualexgeneral} can be generalized to show that it is the alternating
product of the orders of the $\kt$--torsion submodules of $H_*(X;\kt^d)$ (cf. also \cite{KL99}
in the commutative case).
\end{remark}

\section{3--manifolds and 2--complexes}   \label{section:3mfd}

We now restrict ourselves to $\phi$--compatible representations since these have a closer
connection to the topology of a space.

\begin{lemma} \label{lem:alex03}
Let $X$ be a connected CW--complex with finitely many cells in dimensions zero and one. Let
$\phi \in H^1(X;\Z)$ non--trivial and let $\a:\pi_1(X)\to \gl(\kat,d)$ be a $\phi$--compatible
representation. Then $\Delta^\a_0(t)\ne 0$. If $X$ is in fact an $k$--manifold, then
$\Delta^\a_{k}(t)=1$.
\end{lemma}

We need the following notation. If $A=(a_{ij})$ is an $r\times s$--matrix over $\Z[\pi_1(X)]$
and  $\a:\pi_1(X)\to \gl(R,d)$ a representation. Then we denote by $\a(A)$ the $rd\times
sd$--matrix over $R$ obtained by replacing each entry $a_{ij}\in \Z[\pi_1(X)]$ of $A$ by the
$d\times d$--matrix $\a(a_{ij})$.

\begin{proof}
First equip $X$ with a CW--structure with one $0$--cell and $n$ $1$--cells $g_1,\dots,g_n$. We
denote the corresponding elements in $\pi_1(X)$ by $g_1,\dots,g_n$ as well. Since $\phi$ is
non--trivial there exists at least one $i$ such that $\phi(g_i)\ne 0$. Write
$C_*=C_*(\ti{X})\otimes_{\Z[\pi_1(X)]}\kt^d$. The boundary map $\partial_1:C_1\to C_0$ is
represented by the matrix
\[ (\a(1-g_1),\dots,\a(1-g_n))=(\id-\a(g_1),\dots,\id-\a(g_n)).\]
Since $\a$ is $\phi$--compatible it follows that $\a(1-g_i)=\id-At^{\phi(g_i)}$ for some matrix
$A\in \gl(\k,d)$. The first statement of the lemma now follows from Lemma \ref{lem:detaplusb}.

If $X$ is a closed $k$--manifold then equip $X$ with a CW--structure with one $k$--cell. Since
$\phi$ is primitive and $\phi$--compatible an argument as above shows that $\partial_k:C_k\to
C_{k-1}$ has full rank, i.e. $H_k(X;\kt^d)=0$. Hence $\Delta^\a_{k}(t)=1$. If $X$ is a
$k$--manifold with boundary, then it is homotopy equivalent to a $k-1$--complex, and hence
$H_k(X;\kt^d)=0$.
\end{proof}

\begin{lemma} \label{lem:detaplusb}
Let $\kt$ be a skew Laurent polynomial ring and let $A,B$ be invertible $d\times d$--matrices
over $\K$ and $r \ne 0$. Then $\deg(\det(A+Bt^r))=dr$. In particular $A+Bt^r$ is invertible over
$\ktfield$.
\end{lemma}

We point out that Harvey \cite[Proposition~9.1]{Ha05} proves a related result.

\begin{proof}
We can clearly assume that $r>0$. Let $\{e_1,\dots,e_d\}$ be a basis for $\K^d$. Consider the
projection map $p:\kt^d\to P=\kt^d/(A+Bt^r)\kt^d$. Note that by Lemma \ref{lem:finite} we are
done once we show that $p(e_it^j), i \in \{1,\dots,d\}, j\in \{0,\dots,r-1\}$ form a basis for
$P$ as a right $\K$--vector space.

It follows easily from the fact that $A,B$ are invertible that this is indeed a generating set.
Let $v\in \kt^d\sm \{0\}$. We can write $v=\sum_{i=n}^m v_it^i, v_i\in \K^d$ with $v_n \ne 0,
v_m \ne 0$. Since $A,B$ are invertible it follows that $(A+Bt^r)v$ has terms with $t$--exponent
$n$ and terms with $t$--exponent $m+r$. This observation can be used to show that the above
vectors are linearly independent in $P$.
\end{proof}


\noindent We can now  give the proof of Theorem \ref{mainthm1}.

\begin{proof}[Proof of Theorem \ref{mainthm1}]
 Now let $M$ be a 3--manifold whose boundary is empty
or consists of tori. Note that a standard duality argument shows that $2\chi(M)=\chi(\partial
M)=0$. Let $\phi \in H^1(M;\Z)$ be non--trivial, and $\a:\pi_1(M)\to \gl(\kat,d)$ a
$\phi$--compatible representation.

 We first show that $H_*(M;\ktfield^d)=0$ if and only if
$\Delta_{1}^{\a}(t)\ne 0$. Recall that in  Section \label{ref:proof11} we showed that
$H_i(M;\ktfield^d)=0$ if and only if $\Delta_{i}^{\a}(t)\ne 0$. It now follows from Lemma
\ref{lem:alex03} that $H_i(M;\ktfield^d)=0$ for $i=0,3$. If $\Delta_{1}^{\a}(t)\ne 0$, then
$H_1(M;\ktfield^d)=0$. Since $\chi(H_i(M;\ktfield^d))=d \chi(M)=0$ it follows that
$H_2(M;\ktfield^d)=0$.

 \begin{claim}
\[ \tnphi \geq
\frac{1}{d}\left(-\dim_{\K}\left(H_0^\a(M;\kkat)\right)+\dim_{\k}\left(H_1^\a(M;\kkat)\right)-
\dim_{\K}\left(H_2^\a(M;\kkat\right)\right).
\]
Furthermore this inequality becomes an equality if $(M,\phi)$ fibers over $S^1$ and if $M\ne
S^1\times D^2, M\ne S^1\times S^2$.
\end{claim}

 First note that if $\phi$ vanishes on $X\subset M$ then $\a$ restricted to $\pi_1(X)$ lies in
$\glkk \subset \glkatk$ since  $\a$ is $\phi$--compatible. Therefore $H_i^\a(X;\kkat)\cong
H_i^\a(X;\kk)\otimes_\k \kat$. The proofs of \cite[Theorem~3.1]{FK05} and \cite[Theorem~6.1]{FK05}
 can now easily be translated to this
non--commutative setting. This proves the claim.

\noindent Combining the results of the claim with  Lemma \ref{lem:finite} and Corollary
\ref{cor:taudelta}  we immediately get a proof for Theorem \ref{mainthm1}.
\end{proof}

In order to relate Theorem \ref{mainthm1} to the results of  \cite{Co04,Ha05,Tu02b} we need the
following computations for one--dimensional $\phi$--compatible representations.  Recall that
$\phi \in H^1(X;\Z)$ is called primitive if the corresponding map $\phi:H_1(X;\Z)\to \Z$ is
surjective.

\begin{lemma} \label{lem:alex03-onediml}
Let $X$ be a connected CW--complex with finitely many cells in dimensions zero and one. Let
$\phi \in H^1(X;\Z)$ primitive. Let $\a:\pi_1(X)\to \gl(\kat,1)$ be a $\phi$--compatible
one--dimensional representation. If $\a(\pi_1(X))\subset \gl(\kat,1)$ is cyclic, then
$\Delta^\a_{0}(t)=at-1$ for some $a\in \k\sm \{0\}$. Otherwise $\Delta^\a_{0}(t)=1$.
\end{lemma}

\begin{proof}
Equip $X$ with a CW--structure with one $0$--cell and then consider the chain complex for $X$ as in
Lemma \ref{lem:alex03}. The lemma now follows easily from the observation that in $\kt$ we have
$\gcd(1-at,1-bt)=1$ if $a\ne b\in \k$.
\end{proof}


\begin{lemma} \label{lem:alex2-onediml}
Let  $X$ be a 3--manifold with empty or toroidal boundary or let $X$ be a 2--complex with
$\chi(X)=0$. Let $\phi \in H^1(M;\Z)$ non--trivial. Let $\a:\pi_1(M)\to \gl(\kat,1)$ be a
$\phi$--compatible one--dimensional representation. Assume that $\Delta^\a_{1}(t)\ne 0$. If $X$ is
a closed 3--manifold, then $\deg(\Delta^\a_{2}(t))=\deg(\Delta^\a_{0}(t))$, otherwise
$\Delta^\a_{2}(t)=1$.
\end{lemma}

\begin{proof}
First assume that $X$ is a 3--manifold. Then the lemma follows from combining \cite[Sections 4.3
and 4.4]{Tu02b} with \cite[Lemmas 4.7 and 4.9]{FK05}. Note that the results of \cite{FK05} also
hold in the non--commutative setting. If $X$ is a 2--complex then the argument in the proof of
Theorem \ref{mainthm1} shows that $H_2(X;\ktfield^d)=0$. But since $X$ is a 2--complex we have
$H_2(X;\kt^d)\subset H_2(X;\ktfield^d)$, hence $H_2(X;\kt^d)=0$ and $\Delta^\a_2(t)=1$.
\end{proof}

\begin{remark}
Note that we can not  apply the duality results of \cite[Lemma~4.12 and Proposition~4.13]{FK05}
since the natural involution on $\Z[G]$ does not necessarily extend to an involution on $\kt$, i.e.
the representation $\Z[G]\to \kt$ is not necessarily unitary.
\end{remark}

 It now follows immediately from Lemma
\ref{lem:alex03-onediml} and \ref{lem:alex2-onediml} and the discussion in Section
\ref{section:ex} that Theorem \ref{mainthm1} contains the results of \cite{Mc02}, \cite{Co04},
\cite{Ha05}, \cite{Tu02b} and \cite{FK05}.

\begin{remark}
Given a 2--complex $X$ Turaev  \cite{Tu02a} defined a norm $||-||_X:H^1(X;\R)\to \R$,  modelled
on the definition of the Thurston norm of a 3--manifold. In \cite{Tu02a} and \cite{Tu02b} Turaev
gives lower bounds for the Turaev norm which have the same form as certain lower bounds for the
Thurston norm. Going through the proofs in \cite{FK05}  it is not hard to see that the obvious
version of Theorem \ref{mainthm1} for 2--complexes also holds.

If $M$ is a 3--manifold with boundary, then it is homotopy equivalent to a 2--complex $X$. It is
not known whether the Thurston norm of $M$ agrees with the Turaev norm on $X$. But the fact that
Theorem \ref{mainthm1} holds in both cases, and the observation that $\deg(\tau(X,\a))$ is a
homotopy invariant by Theorem \ref{thm:taualexgeneral} suggests that they do in fact agree.
\end{remark}

\section{Examples for skew fields and $\phi$--compatible representations} \label{section:ex}

\subsection{Skew fields of group rings} \label{section:skewfield}
A group $G$ is called locally indicable if for every finitely generated subgroup $U\subset G$
there exists a non--trivial homomorphism $U\to \Z$.

\begin{theorem}\label{thm:tfa}
Let $G$ be a locally indicable and amenable group and let $\F$ be a commutative field. Then the
following hold.
\bn
\item  $\F[G]$  is an Ore domain, in particular it embeds in its classical right ring of
quotients $\K(G)$.
\item $\K(G)$ is flat over $\F[G]$.
\en
\end{theorem}

It follows from \cite{Hi40}  that $\F[G]$ has no zero divisors. The first part now follows from
\cite{Ta57} or \cite[Corollary~6.3]{DLMSY03}. The second part is a well--known property of Ore
localizations (cf. e.g. \cite[p.~99]{Ra98}). We call $ \K(G)$ the \emph{Ore localization} of
$\F[G]$.

A group $G$ is called poly--torsion--free--abelian (PTFA) if there exists a filtration
\[ 1=G_0 \subset G_1\subset \dots \subset G_{n-1}\subset G_n=G \]
such that $G_{i}/G_{i-1}$ is torsion free abelian. It is well--known that PTFA groups are
amenable and locally indicable (cf. \cite{St74}). The group rings of PTFA groups played an
important role in \cite{COT03}, \cite{Co04} and \cite{Ha05}.

\subsection{Examples for $\phi$--compatible representations}
\label{section:examples} Let $X$ be a connected CW--complex and $\phi \in H^1(X;\Z)$. We give
examples of $\phi$--compatible representations.

Let $\F$ be a commutative field. Note that $\phi \in H^1(X;\Z) \cong \hom(H_1(X;\Z),\langle
t\rangle)$ induces a $\phi$--compatible representation $\phi:\Z[\pi_1(X)]\to \ft$. Furthermore
if $\a:\pi_1(X) \to \gl(\f,d)$ is a representation, then $\pi_1(X)$ acts via $\a\otimes \phi$ on
the $\ft$--module $\fd\otimes_\f \ft\cong \ft^d$. We therefore get a representation
 $\a \otimes \phi:\pi_1(X) \to \gl(\ft,d)$, which is
clearly $\phi$--compatible. In this particular case Theorem \ref{mainthm1} was proved in
\cite{FK05}.

To describe the $\phi$--compatible representations of Cochran \cite{Co04} and Harvey
\cite{Ha05,Ha06} we need the following definition.

\begin{definition} Let $\pi$ be a group, $\phi:\pi \to \Z$ an epimorphism and
$\varphi:\pi\to G$ an epimorphism to a locally indicable and amenable group $G$ such that there
exists a map $\phi_G:G\to \Z$ (which is necessarily unique) such that
 \[   \xymatrix {
 \pi \ar[dr]_{\phi} \ar[r]^{\varphi} &G \ar[d]^{\phi_G}\\& \Z }
     \]
 commutes. Following \cite[Definition~1.4]{Ha06} we call $(\varphi,\phi)$
an {\em admissible pair}. If $\phi_G$ is an isomorphism, then $(\varphi,\phi)$ is called {\em
initial}.
\end{definition}

Now let $(\varphi:\pi_1(X)\to G,\phi)$ be an admissible pair for $\pi_1(X)$. In the following we
denote $\ker\{\phi:G\to \Z\}$ by $G'(\phi)$. When the homomorphism $\phi$ is understood we will
write $G'$ for $G'(\phi)$. Clearly $G'$ is still a locally indicable and amenable group. Let
$\F$ be any commutative field and $\K(G')$ the Ore localization of $\F[G]$. Pick an element $\mu
\in G$ such that $\phi(\mu)=1$. Let $\g:\K(G')\to \K(G')$ be the homomorphism given by
$\g(a)=\mu a\mu^{-1}$. Then we get a representation
\[ \ba{rcl} G&\to& \gl(\K(G')_\g\tpm,1)\\[1mm]
    g &\mapsto& (g\mu^{-\phi(g)}t^{\phi(g)}).\ea \]
    It is clear that $\a:\pi_1(X)\to G\to
\gl(\K(G')_{\g}\tpm,1)$ is $\phi$--compatible. Note that the ring $\K(G')_\g\tpm$ and hence the
above representation depends on the choice of $\mu$. We will nonetheless suppress $\g$ in the
notation since different choices of splittings give isomorphic rings. We often make use of the
fact that $f(t)g(t)^{-1}\to f(\mu)g(\mu)^{-1}$ defines an isomorphism $\K(G')(t) \to \K(G)$ (cf.
\cite[Proposition~4.5]{Ha05}). Similarly $\Z[G']\tpm \xrightarrow{\cong} \Z[G]$.

An important example of admissible pairs is provided by Harvey's rational derived series of a group
$G$ (cf. \cite[Section~3]{Ha05}). Let $G_r^{(0)}=G$ and define inductively
\[ G_r^{(n)}=\big\{ g\in G_r^{(n-1)} | \, g^k \in \big[G_r^{(n-1)},G_r^{(n-1)}\big] \mbox{ for some }k\in \Z \sm \{0\} \big\}.\]
Note that $G_r^{(n-1)}/G_r^{(n)}\cong
\big(G_r^{(n-1)}/\big[G_r^{(n-1)},G_r^{(n-1)}\big]\big)/\mbox{$\Z$--torsion}$.
 By \cite[Corollary~3.6]{Ha05} the quotients $G/G_r^{(n)}$ are PTFA groups
for any $G$ and any $n$. If $\phi:G\to \Z$ is an epimorphism, then  $(G\to G/G_r^{(n)},\phi)$ is
an admissible pair for $(G,\phi)$ for any $n>0$.

For example if $K$ is a knot, $G=\pi_1(S^3\sm K)$, then it follows from \cite{St74} that
$G^{(n)}_r=G^{(n)}$, i.e. the rational derived series equals the ordinary derived series (cf. also
\cite{Co04} and \cite{Ha05}).

\begin{remark}
Recall that for a knot $K$ the knot exterior $S^3\sm \nu K$ is denoted by $X(K)$. Let
$\pi=\pi_1(X(K))$ and let $\phi\in H^1(X(K);\Z)$ primitive. Then
\[\ol{\delta}_n(K)=\dim_{\K(\pi'/(\pi')_r^{(n)})}(H_1(X(K),\K(\pi'/(\pi')_r^{(n)})\tpm)\]
is a knot invariant for $n>0$. Cochran \cite[p.~395, Question~5]{Co04} asked whether $K \mapsto
\ol{\delta}_n(K)$ is of finite type.

Eisermann  \cite[Lemma~7]{Ei00} shows that the genus is not a finite type knot invariant.
Cochran \cite{Co04} showed that $\ol{\delta}_n(K)\leq 2\, \genus(K)$ (cf. also Theorem
\ref{mainthm1} together with Corollary \ref{cor:taudelta} and Lemmas \ref{lem:alex03-onediml}
and \ref{lem:alex2-onediml}). Eisermann's argument can now be used to show that $K\mapsto
\ol{\delta}_n(K)$ is not of finite type either.
\end{remark}


Let $X$ be again be a  connected CW--complex and $\phi \in H^1(X;\Z)$.  The two types of
$\phi$--compatible representations given above can be combined as follows. Let $\a:\pi_1(X)\to
\gl(\F,d)$ be a representation and let $\varphi:\pi_1(X)\to G$ be a homomorphism such that
$(\varphi,\phi)$ is an admissible pair. Denote the Ore localization of $\F[G']$ by $\K(G')$.
Then $\pi_1(X)$ acts via $\varphi \otimes\a $ on $\K(G')\tpm \otimes_{\f}\fd \cong
\K(G')\tpm^d$. We therefore get a $\phi$--compatible representation $ \varphi\otimes
\a:\pi_1(X)\to \gl(\K(G')\tpm,d)$.

\section{Comparing different $\phi$--compatible maps}\label{section:comparison}

\noindent We now recall a definition from \cite{Ha06}.

\begin{definition}
Let $\pi$ be a group and $\phi:\pi\to \Z$ an epimorphism. Furthermore let $\varphi_1:\pi \to
G_1$ and $\varphi_2:\pi \to G_2$ be epimorphisms to  locally indicable and amenable groups $G_1$
and $G_2$. We call $(\varphi_1,\varphi_2,\phi)$ an {\em admissible triple} for $\pi$ if there
exist epimorphisms $\varphi_{2}^1:G_1\to G_2$ (which is not an isomorphism) and $\phi_2:G_2\to
\Z$ such that $\varphi_2=\varphi_{2}^1\circ \varphi_1$, and $\phi=\phi_2\circ \varphi_2$.
\end{definition}

\noindent The situation can be summarized in the following diagram
 \[   \xymatrix { &G_1\ar[d]^{\varphi_2^1}\\ \pi \ar[dr]_{\phi} \ar[ur]^{\varphi_1} \ar[r]^{\varphi_2}
&G_2 \ar[d]^{\phi_2}\\& \Z. }
     \]
Note that in particular $(\varphi_i,\phi), i=1,2$ are admissible pairs for $\pi$. Given an
admissible triple we can pick splittings $\Z\to G_i$ of  $\varphi_i, i=1,2$ which make the
following diagram commute:
 \[   \xymatrix { \Z \ar[r] \ar[dr]& G_1 \ar[d]^{\varphi_2^1} \\&G_2.  }    \]
We therefore get an induced commutative diagram of
 ring homomorphisms
 \[   \xymatrix { \Z[\pi] \ar[r] \ar[dr]& \Z[G_1']\tpm \ar[d]^{ \varphi_{2}^1} \\& \Z[G_2']\tpm. }    \]
 Note that we are suppressing the notation for the
 twisting in the skew Laurent polynomial rings.
Denote the $\phi$--compatible maps $\Z[\pi]\to \K(G_i')\tpm, i=1,2$ by $\varphi_i$ as well. For
convenience we recall Theorem \ref{mainthm2}.

\newtheorem*{thm1}{Theorem~\ref{mainthm2}}
 \begin{thm1} \it{
Let $M$ be a 3--manifold whose boundary is a (possibly empty) collection of tori or let $M$ be a
2--complex with $\chi(M)=0$. Let $\a:\pi_1(M)\to \gl(\F,d)$ be a representation and
$(\varphi_1,\varphi_2,\phi)$ an admissible triple for $\pi_1(M)$. If $\tau(M,\varphi_2\otimes
\a)\ne 0$, then $\tau(M,\varphi_1\otimes \a)\ne 0$. Furthermore in that case
\[ \deg(\tau(M,\varphi_1\otimes \a))\geq  \deg(\tau(M,\varphi_2\otimes \a)).\]
}
\end{thm1}


\subsection{Proof of Theorem \ref{mainthm2} for closed 3--manifolds} \label{section:proof3d}
In this section let $M$ be a closed 3--manifold. Choose a triangulation of $M$. Let $T$ be a
maximal tree in the $1$-skeleton of the triangulation and let $T'$ be a maximal tree in the dual
1-skeleton. Following \cite[Section 5]{Mc02} we collapse $T$ to form a single 0-cell and join
the 3-simplices along $T'$ to form a single 3-cell. Since $\chi(M)=0$ the number $n$ of 1--cells
equals the number of 2--cells. Consider the chain complex of the universal cover $\ti{M}$:
\[ 0\to C_3(\ti{M})^1\xrightarrow{\partial_3} C_2(\ti{M})^n\xrightarrow{\partial_2}
C_1(\ti{M})^n\xrightarrow{\partial_1} C_0(\ti{M})^1\to 0,\]
 where the supscript indicates the
rank over $\Z[\pi_1(M)]$. Picking appropriate lifts of the (oriented) cells of $M$ to cells of
$\ti{M}$ we get bases $\ti{\s}_i=\{\ti{\s}^1_i,\dots,\ti{\s}^{r_i}_i\}$ for the
$\Z[\pi_1(M)]$--modules $C_i(\ti{M})$, such that if $A_i$ denotes the matrix corresponding to
$\partial_i$, then $A_1$ and $A_3$ are of the form
\[ \ba{rcl} A_3&=&(1-g_1,\dots,1-g_n)^t, \quad g_i \in \pi_1(M)\\
  A_1&=&(1-h_1,\dots,1-h_n), \quad h_i \in \pi_1(M). \ea \]
Clearly $\{h_1,\dots,h_n\}$ is a generating set for $\pi_1(M)$. Since $M$ is a closed
3--manifold $\{g_1,\dots,g_n\}$ is a generating set for $\pi_1(M)$ as well. In particular we can
find
  $k,l\in \{1,\dots,n\}$ such that $\phi(g_k)\ne 0, \phi(h_l)\ne 0$.

In the following we write $\a_i=\varphi_i \otimes \a:\pi_1(M)\to \gl(\K(G_i')\tpm \otimes
\F^d)\to \gl(\K(G_i')\tpm,d)$, $i=1,2$ and we write $\varphi=\varphi^1_2$.

\begin{lemma} \label{lem:b1b3}
\[ \ba{rcccl} \\deg(\a_1(1-h_l))&=&\deg(\a_2(1-h_l))&=&d |\phi(h_l)| \\
 \deg(\a_1(1-g_k))&=&\deg(\a_2(1-g_k))&=&d |\phi(g_k)| .\ea \]
In particular the matrices $\a_i(1-h_l), \a_i(1-g_k)$ are invertible over $\K(G_i')(t)$ for
$i=1,2$.
\end{lemma}

\begin{proof}
Note that $\a_i(1-h_l)=\id -\a_i(h_l)$, $\a_i(1-g_k)=\id-\a_i(g_k)$ and that $\phi(h_l)\ne 0$,
$\phi(g_k)\ne 0$. The lemma now follows from Lemma  \ref{lem:detaplusb} since $\a_1$ and $\a_2$
are $\phi$--compatible.
\end{proof}

\noindent Denote by  $B$ the result of deleting the $k$--th column and the $l$--row of $A_2$.

\begin{lemma} \label{lem:comptau}
$\tau(M,\a_i)\ne 0$ if and only if $\a_i(B)$ is invertible. Furthermore if $\tau(M,\a_i)\ne 0$,
then
\[ \tau(M,\a_i)= \a_i(1-g_k)^{-1}\a_i(B)\a_i(1-h_l)^{-1}\in {K}_1(\K(G_i')(t))/\pm \a_i(\pi_1(M)). \]
\end{lemma}

\begin{proof}
Denote the standard basis of $\K(G_i')(t)^d$ by $e_1,\dots,e_d$. We equip
$C_j=C^{\a_i}_j(M;\K(G_i')(t)^d)$ with  the ordered bases
 \[ \mathcal{C}_j=\{\ti{\s}^1_j\otimes e_1,\dots,\ti{\s}^1_j\otimes e_d,\dots,
 \ti{\s}_i^{r_j}\otimes e_1,\dots,\ti{\s}_i^{r_j}\otimes e_d\}.\]
Now let \[ \ba{rcl} \xi_0&=&\emptyset,\\
\xi_1&=&\{ld+1,\dots,l(d+1)\},\\
\xi_2&=&\{1,\dots,nd\}\sm \{kd+1,\dots,(k+1)d\},\\
\xi_3&=&\{1,\dots,d\}.\ea \] Then $\xi=(\xi_0,\xi_1,\xi_2,\xi_3)$ is a $\tau$--chain for for
$C_*$. Note that $A_1(\xi)=\a_i(1-h_l), A_2(\xi)=\a_i(B)$ and $A_3(\xi)=\a_i(1-g_k)$. Clearly
$A_1(\xi)$ and $A_3(\xi)$  are invertible by Lemma \ref{lem:b1b3}. The proposition now follows
immediately from Theorem \ref{thm:comptauturaev}.
\end{proof}

Now assume that $\tau(M,\a_2)\ne 0$. Then $\a_2(B)$ is invertible over $\K(G_2')(t)$ by Lemma
\ref{lem:comptau}. Note that $\a_i(B)$ is defined over $\Z[G_i']\tpm\subset \K(G_i')(t)$. In
particular $\a_2(B)=\varphi(\a_1(B))$. It follows from the following lemma that $\a_1(B)$ is
invertible as well.

\begin{lemma} \label{lem:injinv}
Let $P$ be an $r\times s$--matrix over $\Z[G_1']\tpm$. If
\[ \ba{rcl} \Z[G_2]^s &\to& \Z[G_2]^r\\
v&\mapsto &\varphi(P)v \ea \] is invertible (injective) over $\K(G_2')(t)$, then $P$ is
invertible (injective) over $\K(G_1')(t)$.
\end{lemma}

\begin{proof}
Assume that multiplication by $\varphi(P)$ is injective over $\K(G_2')(t)$. Since $\Z[G_2]\to
\K(G_2')(t)=\K(G_2)$ is injective it follows that $\varphi(P):\Z[G_2]^s\to \Z[G_2]^r$ is
injective. By Proposition \ref{prop:strebel} the map $P:\Z[G_1]^s\to \Z[G_1]^r$ is injective as
well. Since $\K(G_1')(t)=\K(G_1)$ is flat over $\Z[G_1]$ it follows that $P:\K(G_1')(t)^s\to
\K(G_1')(t)^r$ is injective.

If $\varphi(P)$ is invertible over the skew field $\K(G_2')(t)$, then $r=s$. But an injective
homomorphism between vector spaces of the same dimension over a skew field is in fact an
isomorphism. This shows that $P$ is invertible over $\K(G_1')(t)$.
\end{proof}

\begin{proposition}\label{prop:strebel}
If $G_1$ is locally indicable, and if $\Z[G_1]^s\to \Z[G_1]^r$ is a map  such that
$\Z[G_1]^s\otimes_{\Z[G_1]}\Z[G_2]\to \Z[G_1]^r\otimes_{\Z[G_1]}\Z[G_2]$ is injective, then
$\Z[G_1]^s\to \Z[G_1]^r$ is injective as well.
\end{proposition}

\begin{proof}
Let $K=\ker\{\varphi:G_1\to G_2\}$. Clearly $K$ is again locally indicable. Note that $\Z[G_1]^s\to
\Z[G_1]^r$ can also be viewed as a map between free $\Z[K]$--modules. Pick any right inverse
$\l:G_2\to G_1$ of $\varphi$. It is easy to see that  $g\otimes h \mapsto g\l(h)\otimes 1, g\in
G_1, h\in G_2$ induces an isomorphism
\[  \Z[G_1]\otimes_{\Z[G_1]} \Z[G_2]\to \Z[G_1]\otimes_{\Z[K]}\Z. \]
By assumption $\Z[G_1]^s\otimes_{\Z[K]} \Z \to \Z[G_2]^r\otimes_{\Z[K]} \Z$ is injective. Since
$K$ is locally indicable it  follows immediately from \cite{Ge83} or \cite{HS83} (cf. also
\cite{St74} for the case of PTFA groups) that $\Z[G_1]^s\to \Z[G_1]^r$ is injective.
\end{proof}

By Lemma \ref{lem:comptau} we now showed that if $\tau(M,\a_2)\ne 0$, then $\tau(M,\a_1)\ne 0$.
Furthermore
\[  \deg(\tau(M,\a_i))=\deg(\a_i(B))
-\deg(\a_i(1-g_k))-\deg(\a_i(1-h_l)), i=1,2.\]
 Theorem \ref{thm:monoton} now follows immediately
from Lemma \ref{lem:b1b3} and from the following proposition.

\begin{proposition} \label{prop:functorial}
 Let $P$ be an $r\times r$--matrix over $\Z[G_1']\tpm$. If $\varphi(P)$ is
invertible then
\[ \deg(P)\geq \deg(\varphi(P)).\]
\end{proposition}

\begin{remark}
\bn
\item
If $\varphi:R\to S$ is a homomorphism of \emph{commutative} rings, and if $P$ is a matrix over
$R\tpm$, then clearly
\[ \ba{rcl} \deg(P)=\deg(\det(P))&\geq& \deg(\varphi(\det(P)))=\deg(\det(\varphi(P)))\\
&=&\deg(\varphi(P)).\ea \] Similarly, several other results in this paper, e.g. Theorem
\ref{thm:welldef} and Lemma \ref{lem:b1b3} are clear in the commutative world, but require more
effort in our noncommutative setting.
\item If
\[ (\Z[G_1'], \{ f \in \Z[G_1'] | \varphi(f)\ne 0 \in \Z[G_2'] \}) \]
has the Ore property, then one can give an elementary proof of the proposition by first
diagonalizing over $\K(G_2')$ and then over $\K(G_1')$. Since this is not known to be the case,
we have to give a more indirect proof.
\en
\end{remark}

\noindent The following proof is based on arguments in \cite{Co04} and \cite{Ha06}.

\begin{proof}[Proof of Proposition  \ref{prop:functorial}]
Let $s=\deg(\varphi(P))$. Pick a map $f:\Z[G_1']^s\to \Z[G_1']\tpm^r$ such that the induced map
\[ \K(G_2')^s \to \K(G_2')\tpm^r \to \K(G_2')\tpm^r/\varphi(P)\K(G_2')\tpm^r \]
is an isomorphism. Denote by $0\to C_1\xrightarrow{P}C_0\to0$ the complex
\[ 0\to\Z[G_1']\tpm^r \xrightarrow{P} \Z[G_1']\tpm^r \to 0,\]
and denote by $0\to D_0\to 0$ the complex with $D_0=\Z[G_1']^s$. We have a chain map $D_*\to
C_*$ given by $f:D_0\to C_0$. Denote by $\mbox{Cyl}(D_*\xrightarrow{f}C_*)$ the mapping cylinder
of the complexes. We then get a short exact sequence of complexes
\[ 0 \to D_*\to \mbox{Cyl}(D_*\xrightarrow{f}C_*)\to \mbox{Cyl}(D_*\xrightarrow{f}C_*)/D_*\to 0.\]
More explicitly we get the following commutative diagram:
\[ \xymatrix{
&0\ar[r]\ar[d]&C_1\oplus D_0 \ar[r]^{{\tiny \bp \id & \id\ep}}
 \ar[d]
 ^{{\tiny{\bp P&-f\\0&\id\ep}}} & C_1\oplus D_0\ar[r] \ar[d]
 ^{{\tiny ( P \,-f)}}
&0\ar[d]\\
0\ar[r]&D_0 \ar[r]^{\tiny{\bp 0 & \id \ep}} &C_0\oplus D_0 \ar[r]& C_0 \ar[r]&0. }
\]
Recall that $\mbox{Cyl}(D_*\xrightarrow{f}C_*)$ and $C_*$ are chain homotopic. Using the
definition of $f$ we therefore see that
\[ f:H_0(D_*;\K(G_2'))\to H_0(\mbox{Cyl}(D_*\xrightarrow{f}C_*),\K(G_2'))\]
is an isomorphism. Since $P$ is invertible over $\K(G_2')(t)$ it follows that
$H_1(\mbox{Cyl}(D_*\xrightarrow{f}C_*);\K(G_2'))=0$. It follows from the long exact homology
sequence corresponding to the above short exact sequence of chain complexes that
$H_1(\mbox{Cyl}(D_*\xrightarrow{f}C_*)/D_*;\K(G_2'))=0$, i.e. the matrix $\bp P &-f\ep$ is
injective over $\K(G_2')$. It follows from Lemma \ref{lem:injinv} that
$H_1(\mbox{Cyl}(D_*\xrightarrow{f}C_*)/D_*;\K(G_1'))=0$ as well.
 Again looking at the long exact homology sequence we get that
 \[ f:H_0(D_*;\K(G_1'))\to H_0(\mbox{Cyl}(D_*\xrightarrow{f}C_*);\K(G_1'))=H_0(C_*;\K(G_1'))\]
 is an injection. Hence
 \[ \ba{rcl}
 \deg(\varphi(P))=s&=&\dim_{\K(G_2')}(H_0(D_*;\K(G_2')))\\
 &=&\dim_{\K(G_1')}(H_0(D_*;\K(G_1')))\\
& \leq &\dim_{\K(G_1')}(H_0(C_*;\K(G_1')))\\
&=&\deg(P).\ea \]
 \end{proof}

\subsection{Proof of Theorem \ref{mainthm2} for  3--manifolds with boundary and for 2--complexes}

First let $X$ be a finite connected 2--complex with $\chi(X)=0$. We can give $X$ a CW--structure
with one 0--cell. If $n$ denotes  the number $n$ of 1--cells, then $n-1$ equals the number of
2--cells. Now consider the chain complex of the universal cover $\ti{X}$:
\[ 0\to C_2(\ti{X})^{n-1}\xrightarrow{\partial_2}
C_1(\ti{X})^n\xrightarrow{\partial_1} C_0(\ti{X})^1\to 0.\] As in Section \ref{section:proof3d}
we pick lifts of the cells of $X$ to cells of $\ti{X}$ to get bases such that if $A_i$ denotes
the matrix corresponding to $\partial_i$, then
\[ A_1=(1-h_1,\dots,1-h_n).  \]
Note that $\{h_1,\dots,h_n\}$ is a generating set for $\pi_1(X)$. Let $l\in \{1,\dots,n\}$ such
that $\phi(l)\ne 0$. The proof of Lemma \ref{lem:comptau} can easily be modified to prove the
following.

\begin{lemma} \label{lem:comptau2}
Denote by $B$ the result of deleting the $l$--row of $A_2$. Then $\tau(X,\a)\ne 0$ if and only
if $\a(B)$ is invertible. Furthermore if $\tau(X,\a)\ne 0$, then
\[ \tau(X,\a)= \a(B)\a(1-h_l)^{-1}\in K_1(\K(G_i')(t))/\pm \a_i(\pi_1(X)). \]
\end{lemma}

The proof of Theorem \ref{mainthm2} for closed manifolds can now easily be modified to cover the
case of 2--complexes $X$ with $\chi(X)=0$.

Now let $M$ be again a 3--manifold whose boundary consists of a non--empty set of tori. A
duality argument shows that $\chi(M)=\frac{1}{2}\chi(\partial(M))=0$. Clearly $M$ is homotopy
equivalent to a 2--complex. Reidemeister torsion is not a homotopy invariant but the following
lemma still allows us to reduce the case of a 3--manifold with boundary to the case of a
2--complex.

\begin{lemma}\cite[p.~56 and Theorem~9.1]{Tu01} \label{lem:mx}
Let $M$ be a 3--manifold with boundary. Then there exists a 2--complex $X$ and a \emph{simple}
homotopy equivalence $M\to X$. In particular, if $\a:\pi_1(X)\cong \pi_1(M) \to \gl(R,d)$ is a
representation such that $H_*(X,R^d)=0$, then
\[ \tau(M,\a) =\tau(X,\a) \in K_1(R)/\pm \a(\pi_1(M)). \]
\end{lemma}

Theorem \ref{mainthm2} for 3--manifolds with boundary now follows from Theorem \ref{mainthm2}
for 2--complexes $X$ with $\chi(X)=0$.

\section{Harvey's monotonicity theorem for groups}\label{section:monotone} \label{section:Harvey}
Let $\pi$ be a finitely presented group and let $(\varphi:\pi\to G,\phi:\pi\to \Z)$ be an
admissible pair for $\pi$. Consider $G'=G'(\phi_G)$ and pick a splitting $\Z \to G$ of $\phi_G$.
As in Section \ref{section:examples} we can consider the skew Laurent polynomial ring
$\K(G')\tpm$ together with the $\phi$--compatible representation $\pi \to \gl(\K(G')\tpm,1)$.

Following \cite[Definition~1.6]{Ha06} we define
 $\ol{\delta}_{G}(\phi)$ to be zero if
$H_1(\pi,\K(G')\tpm)$ is not $\K(G')\tpm$--torsion and
\[ \ol{\delta}_{G}(\phi)=\dim_{\K(G')}(H_1(\pi,\K(G')\tpm)) \]
otherwise.
We give an alternative proof for the following result of Harvey \cite[Theorem~2.9]{Ha06}.

\begin{theorem} \label{thm:harvey1}
If $\pi=\pi_1(M)$, $M$ a closed 3--manifold, and if $(\varphi_1:\pi\to G_1,\varphi_2:\pi \to
G_2,\phi)$ is an admissible triple for $\pi$, then
\[ \ba{rcll} \ol{\delta}_{G_1}(\phi) &\geq& \ol{\delta}_{G_2}(\phi), &\mbox{ if }(\varphi_2,\phi)
\mbox{ is not initial},\\
\ol{\delta}_{G_1}(\phi) &\geq& \ol{\delta}_{G_2}(\phi)-2, &\mbox{ otherwise. }\ea \]
\end{theorem}

\begin{proof}
We  clearly only have to consider the case that $\ol{\delta}_{G_2}(\phi)>0$.
 We can build $K(\pi,1)$ by adding $i$--handles to $M$ with $i\geq 3$. It
therefore follows that for the admissible pairs  $(\varphi_i:\pi\to G_i,\phi)$ we have
\[  \ol{\delta}_{G_i}(\phi)=\dim_{\K(G_i')}(H_1(K_1(\pi,1);\K(G_i')\tpm))=\dim_{\K(G_i')}(H_1(M;\K(G_i')\tpm)).\]
We combine this equality with Theorem \ref{thm:monoton}, Corollary \ref{cor:taudelta} and Lemmas
\ref{lem:finite}, \ref{lem:alex03-onediml}, \ref{lem:alex2-onediml}. The theorem follows now
immediately from the observation that $\im\{ \pi_1(M)\to G_i\to \gl(\K(G_i')\tpm,1)\}$ is cyclic
if and only if $\phi:G_i\to \Z$ is an isomorphism.
\end{proof}

This monotonicity result gives in particular an obstruction for a group $\pi$ to be the
fundamental group of a closed 3--manifold. For example Harvey \cite[Example~3.2]{Ha06} shows
that as an immediate consequence we get the (well--known) fact that $\Z^m, m \geq 4$ is not a
3--manifold group.

\begin{remark}
In \cite{FK05b} the author and Taehee Kim consider the case  $\pi=\pi_1(M)$, where $M$ is a
closed 3--manifold. Given an admissible pair $(\varphi:\pi\to G,\phi)$ we show (under a mild
assumption) that $\ol{\delta}_{G}(\phi)$ is even, generalizing a result of Turaev
(\cite[p.~141]{Tu86}). Furthermore in \cite{FH06} the author and Shelly Harvey will show that
given $\pi\to G$, $G$ locally indicable and amenable, the map
\[ \ba{rcl} \hom(G,\Z)&\to&\Z\\
 \phi&\mapsto& \ol{\delta}_{G}(\phi)\ea \]
defines a seminorm on $\hom(G,\Z)$.
\end{remark}

Let $\pi$ be a finitely presented group of deficiency at least one, for example $\pi=\pi_1(M)$
where $M$ is a 3--manifold with boundary. Using a presentation of deficiency one we can build a
2--complex $X$ with $\chi(X)=0$ and $\pi_1(X)=\pi$. The same proof as the proof of Theorem
\ref{thm:harvey1} now gives the following theorem of Harvey \cite[Theorem~2.2]{Ha06}. In the
case that $\pi=\pi_1(S^3\sm K)$, $K$ a knot, this was first proved by Cochran \cite{Co04}.

\begin{theorem}
If $\pi$ is a finitely presented group of deficiency one and if $(\varphi_1,\varphi_2,\phi)$ is an
admissible triple for $\pi$, then
\[ \ba{rcll} \ol{\delta}_{G_1}(\phi) &\geq& \ol{\delta}_{G_2}(\phi), &\mbox{ if }(\varphi_2,\phi)
\mbox{ is not initial},\\
\ol{\delta}_{G_1}(\phi) &\geq& \ol{\delta}_{G_2}(\phi)-1, &\mbox{ otherwise. }\ea \]
\end{theorem}
\section{Open questions and problems} \label{section:questions}
 Let $M$ be a 3--manifold and $\phi\in H^1(M;\Z)$.
We propose the following three problems for further study.
\bn
\item
 If
$(\varphi:\pi_1(M)\to G,\phi)$ is an admissible pair for $\pi_1(M)$ and if $\a:\pi_1(M)\to
\gl(\F,d)$ factors through $\varphi$, does it follow that
\[ \frac{1}{d}\deg(\tau(M,\a))\leq
\deg(\tau(M,\Z[\pi_1(M)] \to \K(G')(t)))?\] Put differently, are the Thurston norm bounds of
Cochran and Harvey optimal, i.e. at least as good as the Thurston norm bounds of \cite{FK05} for
any representation factoring through $G$.
\item
It is well--known that in many cases $\deg(\tau(M,\Z[\pi_1(M)] \to \K(G')(t)))<||\phi||_T$ for any
admissible pair $(\varphi:\pi_1(M)\to G,\phi)$. For example this is the case if $K$ is a knot with
$\Delta_K(t)=1$ and $M=X(K)$. It is an interesting question whether invariants can be defined for
any map $\pi_1(M)\to G$, $G$ a (locally indicable) torsion--free group. For example it might be
possible to work with $\U(G)$ the algebra of affiliated operators (cf. e.g. \cite{Re98}) instead of
$\K(G)$. If such an extension is possible, then it is a natural question whether the Thurston norm
is determined by such more general bounds. This might be too optimistic in the general case, but it
could be true in the case of a knot complement.
\item If $(M,\phi)$ fibers over $S^1$, then it is well--known that the corresponding Alexander
polynomial defined over $\zt$ is monic, i.e. the top coefficient is $\pm 1$. Because of the
high--degree of indeterminacy of Alexander polynomials over skew Laurent polynomial rings a
corresponding statement is meaningless. Since Reidemeister torsion has a much smaller
indeterminacy it is potentially possible to use it to extend the above fiberedness obstruction
as in \cite{GKM05}.
%
\en

\end{document}